\documentclass[12pt,a4paper]{amsart}
\usepackage[utf8]{inputenc}	
\usepackage[T1]{fontenc}
\usepackage[in,headings]{fullpage}
\usepackage{amsbsy, amscd, amsfonts, amsmath, amsrefs, amssymb, amsthm}
\usepackage{graphicx}
\usepackage{float}
\usepackage{color}   
\usepackage[colorlinks=true,linkcolor=blue,citecolor=blue,urlcolor=blue]{hyperref}

\newtheorem{theorem}{Theorem}

\newtheorem{lemma}[theorem]{Lemma}

\theoremstyle{definition}

\newtheorem{remark}[theorem]{Remark}

\theoremstyle{remark}

\newcommand{\R}{\mathbf{R}}
\newcommand{\N}{\mathbf{N}}

\renewcommand{\Re}{\mathop{\mathrm{Re}}\nolimits}
\renewcommand{\Im}{\mathop{\mathrm{Im}}\nolimits}
\newcommand{\Rzeta}{\mathop{\mathcal R }\nolimits}

\newfont{\cmbsy}{cmbsy10}
\newfont{\cmmib}{cmmib10}
\newcommand{\Orden}{\mathop{\hbox{\cmbsy O}}\nolimits}
\newcommand{\orden}{\mathop{\hbox{\cmmib o}}\nolimits}

\begin{document}

\title{On the number of zeros of $\Rzeta(s)$.}
\author[Arias de Reyna]{J. Arias de Reyna}
\address{%
Universidad de Sevilla \\ 
Facultad de Matem\'aticas \\ 
c/Tarfia, sn \\ 
41012-Sevilla \\ 
Spain.} 




\email{arias@us.es, ariasdereyna1947@gmail.com}


\begin{abstract}
We prove that the number of zeros $\varrho=\beta+i\gamma$ of $\Rzeta(s)$ with
$0<\gamma\le T$ is given by 
\[N(T)=\frac{T}{4\pi}\log\frac{T}{2\pi}-\frac{T}{4\pi}-\frac12\sqrt{\frac{T}{2\pi}}+\Orden(T^{2/5}\log^2 T).\]
Here $\Rzeta(s)$ is the function that Siegel found in Riemann's papers. Siegel related the zeros of $\Rzeta(s)$ to the zeros of Riemann's zeta function. 
Our result on $N(T)$  improves the result of Siegel. 
\end{abstract}

\maketitle

\section{Introduction}
In his paper on Riemann Nachlass Siegel \cite{Siegel} introduced the auxiliar function (see \cite{A166} for its main properties)
\begin{displaymath}
\Rzeta(s)=\int_{0\swarrow1}\frac{x^{-s} e^{\pi i x^2}}{e^{\pi i x}-
e^{-\pi i x}}\,dx.
\end{displaymath}
He studied the number $N(T)$ of zeros $\varrho=\beta+i\gamma$ of $\Rzeta(s)$ with 
$0<\gamma\le T$ and  proved that $N(T)=\frac{T}{4\pi}\log\frac{T}{2\pi}-\frac{T}{4\pi} +\orden(T)$. If all these zeros of $\Rzeta(s)$ were left of the critical line, the function $\zeta(s)$ will have $\frac{T}{2\pi}\log\frac{T}{2\pi}-\frac{T}{2\pi}$ zeros on the critical line, that is, all but an infinitesimal proportion of zeros of zeta will lie on the critical line.  

After giving the approximate number of zeros to height $T$, Riemann asserted in his paper \cite{R} (see Edwards \cite{Edwards}*{p.~301}) \emph{One finds in fact about this many real roots within these bounds} meaning that almost all nontrivial zeros appear to be in the critical line. Much more explicit is Riemann in a letter to Weierstrass (\cite{R2}*{p.~823--824}), where he said that he thinks Weierstrass will have no problem filling the gaps in his paper (an understatement, it took mathematicians 70 years to fill  these gaps),  except for two Theorems, one of them precisely showing that there are so many zeros in the critical line. Siegel tried  unsuccessfully to use $\Rzeta(s)$ to prove this assertion of Riemann.  If all the zeros of $\Rzeta(s)$ had been to the left of the critical line, the Riemann assertion would follow.   Our computation of zeros \cite{A172} of $\Rzeta(s)$ shows that there are many zeros of $\Rzeta(s)$ to the right of the critical line, perhaps $1/3$ of the zeros of $\Rzeta(s)$ satisfies $\beta>1/2$. 

Here we refine the Siegel result showing that 
\[N(T)=\frac{T}{4\pi}\log\frac{T}{2\pi}-\frac{T}{4\pi}-\frac12\sqrt{\frac{T}{2\pi}}+\Orden(T^{2/5}\log^2 T).\]
The term $\frac12\sqrt{\frac{T}{2\pi}}$ was first guessed in \cite{A172} using the computed zeros of $\Rzeta(s)$, this numerical study suggests that the real error is possibly logarithmic.

To prove our result, we use  the argument principle and my refinements of Siegel asymptotics for $\Rzeta(s)$ (see \cite{A98}, \cite{A100}) and Backlund's lemma. Since I have seen published unnecessarily complicated demonstrations of this lemma and not very neat versions of it, we include here a complete demonstration of this lemma.

\section{Backlund lemma.}

To count the zeros  of holomorphic functions, we have to estimate the real parts of the integrals $\frac{1}{2\pi i}\int f'(z)/f(z)\,dz$. There is an old Lemma of Backlund to obtain this bound. The result follows from the Jensen inequality and can also be substituted by Theorem D in Ingham \cite{I}*{p.~49}. Nevertheless, I have not found a neat version as the one I give of Blacklund's Lemma (compare the statement and proof in \cite{SCH} and \cite{SW}).

\begin{lemma}[Backlund]
Let  $f$ be holomorphic in the disc $|z-a|\le R$. 
Let  $|f(z)|\le M$ for  $|z-a|\le R$. Let  $b$ be a point in the interior of the disc
$0<|b-a|<R$. Assume that  $f$ does not vanish on the segment $[a,b]$, then 
\begin{equation}\label{E:backlund}
\Bigl|\Re\frac{1}{2\pi i}\int_a^b\frac{f'(z)}{f(z)}\,dz\Bigr|\le
\frac{1}{2}\Bigl(\log \frac{M}{|f(a)|}\Bigr)\Bigl(\log\frac{R}{|b-a|}\Bigr)^{-1}.
\end{equation}
\end{lemma}

\begin{proof} 
Since $f(a)\ne0$ and $b-a\ne0$, we have $f(a)=e^{i\alpha}|f(a)|$ and $b-a=e^{i\theta}
|b-a|$, with $\alpha$ and $\theta\in\R$. Note that $0<r:=\frac{|b-a|}{R}<1$, with these notation, define for $|\zeta|\le1$, $g(\zeta)=e^{-i\alpha}f(a+Re^{i\theta}\zeta)$. Then a simple change of variables shows that 
\[\frac{1}{2\pi i}\int_a^b\frac{f'(z)}{f(z)}\,dz=\frac{1}{2\pi i}\int_0^r\frac{g'(\zeta)}{g(\zeta)}\,d\zeta.\]
The problem reduces to the case of $g(\zeta)$ defined in the unit disc and the interval $[0,r]$ with $0<r<1$. We also have that $g(0)=|f(a)|>0$, and $|g(\zeta)|\le M$ for $|\zeta|\le 1$. 
Thus, the problem is reduced to the particular case in which the function  $g(\zeta)$ is holomorphic on the unit disc, $a=0$, $b=r$, $R=1$ and $|g(0)|>0$.

Since  $g$ does not vanish on
$[0,r]$ there is a well-defined holomorphic function  $\log g(z)$ in  a neighborhood of 
this segment. Then we have
\[
I=\Re\frac{1}{2\pi i}\int_0^r\frac{g'(z)}{g(z)}\,dz =\Re\frac{\log g(r)-\log 
g(0)}{2\pi i}= \frac{\arg g(r)-\arg g(0)}{2\pi}.
\]

Introduce now a new function  $h(z)=
\bigl(g(z)+\overline{g(\overline{z})}\bigr)/2$, it is holomorphic on the unit disc, bounded in absolute value by $M$ and does not vanish identically because $h(0)=g(0)>0$. Therefore, $h(z)$ has at most a finite number of zeros in the segment $[0,r]$. Let $(a_j)_{j=1}^n$ be the zeros of $h(z)$ contained in the segment $[0,r]$, repeated according to its multiplicity.  The function $h(x)=\Re(g(x))$ maintains its sign between any two consecutive $a_j$. Its argument can only change in an amount $\pi$ for each of these intervals. Therefore, 
\[I\le \frac{(n+1)\pi}{2\pi}.\]
To estimate $n$, consider the function 
\[H(\zeta):=h(\zeta)\prod_{j=1}^n\frac{1-\overline{a_j}\zeta}{\zeta-a_j}.\]
$H(z)$ is holomorphic on $|\zeta|\le1$ and  is bounded by $M$ on this disc. 
Since all the factors are of modulus $=1$ for $|\zeta|=1$.  By the maximum modulus principle, we have 
$|H(0)|\le M$. This is equivalent to 
\begin{displaymath}
h(0)\Bigl|\prod_{j=1}^n \frac{1}{a_j}\Bigr|\le M.
\end{displaymath}
Since $0<a_j\le r$ for each $j$ and since $h(0)=g(0)$ we obtain
\begin{displaymath}
r^{-n}\le \frac{M}{g(0)}.
\end{displaymath}
so that 
\[n\le \frac{\log \bigl(M/g(0)\bigr)}{\log (1/r)}.\]
Substituting into the above inequality we obtain
\[I\le \frac12+\frac12\frac{\log \bigl(M/g(0)\bigr)}{\log (1/r)}.\]
Note also that $r=|b-a|/R$ in the notation of the statement of the lemma.

To remove the summand $\frac12$, for $m\in\N$ we apply the above reasoning to $f(z)^m$  we obtain
\[m\Bigl|\Re\int_a^b \frac{f'(z)}{f(z)}\,dz\Bigr|\le \frac12+\frac12\Bigl(\log \frac{M^m}{|f(a)|^m}\Bigr)\Bigl(\log\frac{R}{|b-a|}\Bigr)^{-1}\]
Hence,
\[\Bigl|\Re\int_a^b \frac{f'(z)}{f(z)}\,dz\Bigr|\le \frac{1}{2m}+\frac12\Bigl(\log \frac{M}{|f(a)|}\Bigr)\Bigl(\log\frac{R}{|b-a|}\Bigr)^{-1}.\]
Taking the limits for  $m\to+\infty$ we get \eqref{E:backlund}.
\end{proof}

\begin{remark}
The same bound  applies to the integral on  any segment $L$ located in the interior of the circle in a line passing through the center $a$, changing $|a-b|$ by the maximum of $|z-a|$ for $z\in L$. The proof being the same. 
\end{remark}

\section{Applying the argument principle.}

\begin{theorem}
The number $N(T)$ of zeros $\varrho=\alpha+i\beta$ of $\Rzeta(s)$ with $0<\gamma\le T$ satisfies 
\begin{equation}\label{E:main}
N(T)=\frac{T}{4\pi}\log\frac{T}{2\pi}-\frac{T}{4\pi}-\frac12\sqrt{\frac{T}{2\pi}}+
\Orden(T^{2/5}\log^2 T).
\end{equation}
\end{theorem}

\begin{proof}
By Corollary 14 in \cite{A100} we have $|\Rzeta(s)-1|\le \frac34$ for $\sigma\ge2$  and $|s|\ge r_0$. Therefore, there are no zeros of $\Rzeta(s)$ with $\beta\ge2$ except perhaps a finite number of them contained in a closed disc $|s|\le r_0$. By Theorem 11 in \cite{A98} there is a $t_0>0$ and a contant $A>0$ such that 
for $t>t_0$ and $\sigma\le 1-At^{2/5}\log t$ we have
\begin{equation}\label{defU}
\Rzeta(s)=-\chi(s)\eta^{s-1}e^{-\pi i \eta^2}\frac{\sqrt{2}e^{3\pi i/8}\sin\pi\eta}
{2\cos2\pi\eta}(1+U),
\end{equation}
with $|U|<1$ and $\eta=\sqrt{(s-1)/2\pi i}$ the square root determined by $\Re(\eta)+\Im(\eta)>0$. In particular, there are no zeros in this region.

For $0<t\le t_0$ and $\sigma<0$, by Corollary 9 in \cite{A98} the possible zeros  are on the closed disc $|s-1|\le 3528\pi$ or in  $\sigma\ge 1-t\ge 1-t_0$. Therefore, they are finite in number.

The constant $t_0$ of Theorem 11 in \cite{A98} can be taken $t_0>r_0$ and $t_0>3528\pi$, also in such a way that there is no zero of $\Rzeta(s)$ with $\gamma=t_0$ and for later use assume also that $\frac{At^{2/5}\log t}{t}<\frac12$ for $t\ge t_0$.

From all the results quoted, it follows that all the zeros of $\Rzeta(s)$ with $0<\gamma\le t_0$ are contained in a certain rectangle  $[- R, R]\times[0,t_0]$. It follows that $N(t_0)$ is finite. 

Since $N(t_0)$ is finite to prove \eqref{E:main}  we only have to prove that   $N(T)-N(t_0)=\frac{T}{4\pi}\log\frac{T}{2\pi}-\frac{T}{4\pi}-\frac12\sqrt{\frac{T}{2\pi}}+\Orden(T^{2/5}\log^2 T)$.

We want to estimate the number of zeros $\beta+i\gamma$ with $t_0<\gamma\le T$, we can assume that $T$ is chosen so that there is no zero with $\gamma=T$, we also take $t_0$,so that there  is no zero with $\gamma=t_0$. According to the quoted theorems, any  zero with  $t_0<\gamma\le T$ is contained in the set 
\[F=\{s=\sigma +i t\ \colon\ t_0\le t\le T,\quad 1-At^{2/5}\log t\le \sigma\le 2\}.\]
The  closed set $F$ is limited by a contour $L=L_1+L_2+L_3+L_4$ composed of 
$L_1$ the segment that joins $1-At_0^{2/5}\log t_0+it_0$ to $2+it_0$; $L_2$ the segment that joins $2+it_0$ with $2+iT$; $L_3$ the segment that joins $2+iT$ with $1-AT^{2/5}\log T+iT$; and $L_4$ is a curve line given in parametric coordinates by 
$1-At^{2/5}\log t+it$ for $T\ge t\ge t_0$. By the argument principle we have
\[N(T)-N(t_0)=\frac{1}{2\pi i}\int_L\frac{\Rzeta'(s)}{\Rzeta(s)}\,ds.\]
It is sufficient to estimate each of the integrals
\[I_j:=\Re\Bigl\{\frac{1}{2\pi i}\int_{L_j}\frac{\Rzeta'(s)}{\Rzeta(s)}\,ds\Bigr\}.\]
Each of these integrals is equal to the variation of the argument of $\Rzeta(s)$  along the corresponding $L_j$ divided by $2\pi$, which is well defined since $\Rzeta(s)$ do not vanish on the contour of $F$.

(a) \textbf{Estimate of $I_1$ and $I_2$.} 
The variation of the argument along $L_1$ is a constant dependent on $t_0$, which is fixed on this proof. The variation along $L_2$ is at most $\pi$ since for $s\in L_2$
we have $|\Rzeta(s)-1|<3/4$. 

(b) \textbf{Estimate of $I_3$.} 
In segment $L_3$ we apply the Backlund lemma. Take the center of the disc at  $a=2+iT$. Take the circle $D$ with center at $a$ and radius $R=2+2AT^{2/5}\log T$. For $s\in D$, 
$t>16\pi$ and either $\sigma>0$ and then by \cite{A92}*{Prop.~12} we have $|\Rzeta(s)|\le\sqrt{t/2\pi}$ or $\sigma\le0$ and by \cite{A92}*{Prop.~13} we have
\[|\Rzeta(s)|\le \frac{19t}{(2\pi)^{1-\sigma}}\{(1-\sigma)^2+t^2\}^{\frac14-\frac{\sigma}{2}}\le \frac{19t}{\sqrt{2\pi}}\{(\tfrac{1-\sigma}{2\pi})^2+(\tfrac{t}{2\pi})^2\}^{\frac12-\frac{\sigma}{2}}.\]
It follows that the maximum of $|\Rzeta(s)|$ for $s\in D$ is bounded by
\[M=\exp(c T^{2/5}\log^2 T),\]
for some constant $c$. At the center $a$ we have $|\Rzeta(a)|>1/4$. 
Backlund's lemma shows that the variation of the argument of $\Rzeta(s)$ in $L_3$
is  $\ll T^{2/5}\log^2 T$. 

(c) \textbf{Estimate of $I_4$.} 
On the  $L_4$ curve line, the expression \eqref{defU} applies.  The variation of the argument of $\Rzeta(s)$ on this line differs from the variation of the argument of 
\begin{equation}\label{E:factors}
-\chi(s)\eta^{s-1}e^{-\pi i \eta^2}\frac{\sqrt{2}e^{3\pi i/8}\sin\pi\eta}
{2\cos2\pi\eta}
\end{equation}
in a constant less than or equal $\pi$. 

Each factor in \eqref{E:factors} does not vanish in $L_4$. For each $s\in L_4$, we may compute the argument of each factor and sum all to obtain the total argument.
The points of this line are $s=\sigma+it$ where $\sigma=1-At^{2/5}\log t$.
First, we compute $\eta$ and $\log\eta$ with a good approximation in terms of $t$
\begin{equation}\label{etat}
\eta=\Bigl(\frac{t}{2\pi}+\frac{1-\sigma}{2\pi}i\Bigr)^{1/2}=
\sqrt{\frac{t}{2\pi}}\Bigl(1+i\frac{1-\sigma}{2t}+\frac{(1-\sigma)^2}{8t^2}-i\frac{(1-\sigma)^3}{16t^3}-\frac{5(1-\sigma)^4}{128t^4}+\cdots\Bigr)
\end{equation}

\[\log\eta=\frac12\log\frac{t}{2\pi}+\frac{1-\sigma}{2t}i+\frac{(1-\sigma)^2}{4t^2}-i\frac{(1-\sigma)^3}{6t^3}+\cdots\]
with convergent expansions since $(1-\sigma)/t=\frac{At^{2/5}\log t}{t} <\frac12$ for $t\ge t_0$ by the election of $t_0$.
Also, we will need
\begin{align*}
\log s&=\log(\sigma+it)=\tfrac12\log(\sigma^2+t^2)+ i\Bigl(\frac{\pi}{2}+\arctan\frac{|\sigma|}{t}\Bigr)\\=
&\log t+\frac{\sigma^2}{2t^2}-\frac{s^4}{4t^4}+\cdots+i\Bigl(\frac{\pi}{2}-\frac{\sigma}{t}+\frac{\sigma^3}{3t^3}+\cdots\Bigr)
\end{align*}

Estimating $\arg\log\chi(s)$, notice that this is $\Im\log\chi(s)$ and we use the Euler-MacLaurin expansion of $\log\Gamma(s)$ \cite{Edwards}*{p.~109}. 
\begin{align*}
\log\chi(s)&=\log\frac{(2\pi)^s}{2\Gamma(s)\cos(\pi s/2)}=s\log(2\pi)-\log\Gamma(s)-\log(e^{\pi i s/2}+e^{-\pi i s/2})\\
&=s\log(2\pi)-(s-\tfrac12)\log s+s-\tfrac12\log(2\pi)+\Orden(|s|^{-1})-\tfrac{\pi i s}{2}+\Orden(e^{-\pi t}).
\end{align*}
So, with $s=1-At^{2/5}\log t+it$
\begin{multline*}
\Im\log\chi(s)\\ =t\log(2\pi)-(\sigma-\tfrac12)\Bigl(\frac{\pi}{2}-\frac{\sigma}{t}+\frac{\sigma^3}{3t^3}+\cdots\Bigr)-t\Bigl(\log t+\frac{\sigma^2}{2t^2}-\frac{\sigma^4}{4t^4}+\cdots\Bigr)+t+\Orden(t^{-1})-\frac{\pi\sigma}{2}
\end{multline*}
Since $1-\sigma=At^{2/5}\log t$ 
\begin{align*}
\Im\log\chi(s)&=-t\log\frac{t}{2\pi}+t+\frac{\pi}{4}-\pi\sigma
-\frac{\sigma}{2t}+\frac{\sigma^2}{2t}+\Orden(t^{-7/5}\log^3t)\\&=-t\log\frac{t}{2\pi}+t+\Orden(t^{2/5}\log t).\end{align*}
Therefore,
\begin{equation}\label{E:chi}
\arg\chi(s)=-t\log\frac{t}{2\pi}+t+\Orden(t^{2/5}\log t).
\end{equation}

Estimating $\arg \eta^{s-1}$. We have
\begin{align*}
\Im\log\eta^{s-1}&=\Im((s-1)\log\eta)\\
&=\Im\Bigl\{(\sigma-1+it)\Bigl(\frac12\log\frac{t}{2\pi}+\frac{1-\sigma}{2t}i+\frac{(1-\sigma)^2}{4t^2}-i\frac{(1-\sigma)^3}{6t^3}+\cdots\Bigr)\Bigr\}\\
&=-\frac{(1-\sigma)^2}{2t}+\frac{(1-\sigma)^4}{6t^3}+\cdots+\frac{t}{2}\log\frac{t}{2\pi}+\frac{(1-\sigma)^2}{4t}+\cdots\\ &=\frac{t}{2}\log\frac{t}{2\pi}+\Orden(t^{-1/5}\log^2t).
\end{align*}
So, 
\begin{equation}\label{E:etas}
\arg \eta^{s-1}=\frac{t}{2}\log\frac{t}{2\pi}+\Orden(t^{-1/5}\log^2t).
\end{equation}
\begin{equation}\label{E:eta2}
\arg e^{-\pi i\eta^2}=\Im\log e^{-\pi i\eta^2}=\Im\bigl(-\pi i \eta^2\bigr)=\Im((1-s)/2)=-\frac{t}{2}.
\end{equation}

Estimation of the last factor of \eqref{E:factors}
\begin{align*}
\Im\log\Bigl(-e^{3\pi i/8}&\frac{\sin\pi\eta}{\cos2\pi \eta}\Bigr)=-\pi+\frac{3\pi }{8}+\Im\log\frac{e^{\pi i \eta}-e^{-\pi i\eta}}{i(e^{2\pi i \eta}+e^{-2\pi i\eta})}\\
&=-\frac{\pi}{8}+\Im \pi i \eta+\Im\log\frac{1-e^{2\pi i\eta}}{1+e^{4\pi i\eta}}\\
&=-\frac{\pi}{8}+\Im\pi i \sqrt{\frac{t}{2\pi}}\Bigl(1+i\frac{1-\sigma}{2t}+\frac{(1-\sigma)^2}{8t^2}-i\frac{(1-\sigma)^3}{16t^3}+\cdots\Bigr)+\Orden(1)\\
&=\pi\sqrt{\frac{t}{2\pi}}+\Orden(1).
\end{align*}
We have used that $|e^{2\pi i\eta}|<1$ since $\Im(\eta)>0$ although small. Hence 
$\log\frac{e^{\pi i \eta}-e^{-\pi i\eta}}{i(e^{2\pi i \eta}+e^{-2\pi i\eta})}$ can be defined as a difference of two principal logarithms with an imaginary part  with an absolute value less than $\pi/2$, each of them, so that $\Im\log\frac{e^{\pi i \eta}-e^{-\pi i\eta}}{i(e^{2\pi i \eta}+e^{-2\pi i\eta})}$ is in absolute value less than $\pi$.

Combining \eqref{E:chi}, \eqref{E:etas}, \eqref{E:eta2} and the last expression, we obtain
\begin{align*}
\arg&\Bigl(-\chi(s)\eta^{s-1}e^{-\pi i \eta^2}\frac{\sqrt{2}e^{3\pi i/8}\sin\pi\eta}{2\cos2\pi\eta}\Bigr)=
-t\log\frac{t}{2\pi}+t+\Orden(t^{2/5}\log t)\\&
+\frac{t}{2}\log\frac{t}{2\pi}+\Orden(t^{-1/5}\log^2t)
-\frac{t}{2}+\pi\sqrt{\frac{t}{2\pi}}+\Orden(1)\\
&=-\frac{t}{2}\log\frac{t}{2\pi}+\frac{t}{2}+\pi\sqrt{\frac{t}{2\pi}}+
\Orden(t^{2/5}\log^2t).
\end{align*}

It follows that the variation of the argument is 
\[\frac{1}{i}\int_{L_1+L_2+L_3+L_4}\frac{\Rzeta'(s)}{\Rzeta(s)}= 
\frac{T}{2}\log\frac{T}{2\pi}-\frac{T}{2}-\pi\sqrt{\frac{T}{2\pi}}+
\Orden(T^{2/5}\log^2 T)\]
Dividing by $2\pi$ we get \eqref{E:main}.
\end{proof}

\end{document}